\documentclass[11pt]{article}
\usepackage[numbers,sort&compress]{natbib}
\usepackage[colorlinks,
            linkcolor=blue,
            anchorcolor=green,
            citecolor=magenta
            ]{hyperref}
\usepackage{mathrsfs,amssymb,amsfonts}
\usepackage{graphicx}
\usepackage{ifpdf}
\usepackage{amsmath,amssymb,latexsym,color}
\usepackage[mathscr]{eucal}
\usepackage[numbers]{natbib}
\newtheorem{thm}{Theorem}[section]

\newtheorem{prop}[thm]{Proposition}
\newtheorem{lem}[thm]{Lemma}

\newtheorem{cor}[thm]{Corollary}
\newtheorem{example}{Example}
\newtheorem{remark}{Remark}[section]

\newenvironment{eg}{\begin{example} \rm}{\end{example}}

\newcommand{\proof}{{\it Proof.\quad}}
\newcommand{\rc}{{\rm rc}}
\newcommand{\qed}{\hfill\Box\medskip}

\usepackage{CJK}
\begin{document}
\begin{CJK*}{GBK}{song}
\renewcommand{\abovewithdelims}[2]{
\genfrac{[}{]}{0pt}{}{#1}{#2}}

\title{\bf The rainbow connection number of the power graph of a finite group}

\author{Xuanlong Ma$^{\rm a}$\quad Min Feng$^{\rm b, a}$\quad Kaishun Wang$^{\rm a}$\\
{\footnotesize  $^{\rm a}$ \em    Sch. Math. Sci. {\rm \&} Lab. Math. Com. Sys., Beijing Normal University, Beijing, 100875,  China} \\
{\footnotesize  $^{\rm b}$ \em  Department of Applied  Mathematics, Nanjing University of Science and Technology, Nanjing 210094, China }
}
 \date{}
 \maketitle

\begin{abstract}
This paper studies the rainbow connection number of the power graph $\Gamma_G$ of a finite group $G$. We determine the rainbow connection number of $\Gamma_G$ if $G$ has maximal involutions or is nilpotent, and show that the rainbow connection number of $\Gamma_G$ is at most three if $G$ has no maximal involutions. The rainbow connection numbers of  power graphs of some nonnilpotent groups are also given.


\end{abstract}
\footnote{{\em E-mail address:} xuanlma@mail.bnu.edu.cn (X. Ma), fgmn\_1998@163.com (M. Feng), \\
wangks@bnu.edu.cn (K. Wang).}

\section{Introduction}

Given a connected graph $\Gamma$, denote by $V(\Gamma)$ and $E(\Gamma)$ the vertex set and edge set, respectively.
Define a coloring $\zeta: E(\Gamma)\rightarrow \{1,2,\ldots,k\}$, $k\in \mathbb{N}$,
where adjacent edges may be colored the same. A path $P$ is {\em rainbow} if any two edges in $P$ are colored distinct.
If $\Gamma$ has a rainbow
path from $u$ to $v$ for each pair of vertices $u$ and $v$, then $\Gamma$ is {\it rainbow-connected} under the coloring $\zeta$, and $\zeta$ is called a {\it rainbow $k$-coloring} of $\Gamma$.
The {\it rainbow connection number} of $\Gamma$, denoted by rc$(\Gamma)$, is the minimum $k$ for which there exists a rainbow
$k$-coloring of  $\Gamma$.

The rainbow connection number of a graph $\Gamma$ was introduced by Chartrand et al. \cite{CJMZ}. It was showed in \cite{CFMY,LT} that computing $\rc(\Gamma)$ is NP-hard.  Moreover, it has
been proved in \cite{LT}, that for any fixed $t\ge 2$, deciding if $\rc(\Gamma)=t$ is NP-complete.
Some topics on   restrict  graphs are as follows: oriented graphs \cite{DS}, graph products \cite{GMP}, hypergraphs \cite{CRLM}, corona graphs \cite{ES}, line graphs \cite{LSue}, Cayley graphs \cite{LLL},
dense graphs \cite{LlS} and  sparse random graphs \cite{FT}. Most of the results and papers that dealt with it can be found in \cite{LSS}.

In this paper we study the rainbow connection number of the power graph of a finite group.
We always use $G$ to denote a finite group with the identity $e$.
The {\it power graph} $\Gamma_G$ has the vertex set $G$ and two distinct elements are adjacent if one is a power of the other.
The concept of a power graph was introduced in \cite{n1}. Recently, many interesting results on power graphs have been obtained, see \cite{Cam,CGh,CGS,FMW,FMW1,kel21,kel2,kel22,MF}. A detailed list of results and open questions on power graphs can be found in \cite{AKC}. (Since our paper deals only with undirected graphs, for convenience
throughout we use the term ``power graph'' to refer to an undirected power graph defined as above, see also \cite{n1}, Section 3).

A finite group is called a $p$-group if its order is a power of $p$, where $p$ is a prime. In $G$, an element of order $2$ is called an involution.
An involution $x$ is  {\it maximal } if
the only cyclic subgroup containing $x$ is the subgroup $\langle x\rangle$  generated by $x$. Denote by $M_G$ the set of all maximal involutions of $G$. We use $M_G$ to discuss the rainbow connection number of $\Gamma_G$.

This paper is organized as follows.  In Section~\ref{MGn0} we express ${\rm rc}(\Gamma_G)$ in terms of $|M_G|$ if $M_G\neq\emptyset$. In Section~\ref{MG=0}  we show that $\rc(\Gamma_G)\le 3$ if $M_G=\emptyset$. In particular, we determine $\rc(\Gamma_G)$ if $G$ is nilpotent. The rainbow connection numbers of  power graphs of some nonnilpotent groups are also given.

\section{$M_G\ne\emptyset$}\label{MGn0}
In this section we prove the following theorem.

\begin{thm}\label{thm1}
Let $G$ be a finite group of order at least $3$.  Then
$$
{\rm rc}(\Gamma_{G})=\left\{
                                  \begin{array}{ll}
                                    3, & \hbox{if $1 \le |M_G|\le 2$;}\\
                                    |M_G|, & \hbox{if $|M_G|\ge 3$}.
                                  \end{array}
                                \right.
$$
\end{thm}

We begin with the following lemma.

\begin{lem}\label{low}
 ${\rm rc}(\Gamma_{G})\ge |{M}_G|$.
\end{lem}
\proof  Let $M_G=\{z_1,\ldots, z_m\}$.
Observe that $e$ is the unique vertex adjacent to
$z_i$ in $\Gamma_G$, where $i=1,\ldots,m$. Hence, for each pair of maximal involutions $z_i$ and $z_j$,
the path from $z_i$ to $z_j$ is unique, which is $(z_i,e,z_j)$.
Suppose  $\zeta$ is a  rainbow $k$-coloring of $\Gamma_G$. Then
$|\{\zeta(\{z_i,e\}):i=1,\ldots,m\}|=m$, and so $k\ge m$, as desired.
$\qed$

For $x\in G$, let $[x]=\{y\in G: \langle y\rangle=\langle x\rangle\}$. Then $\{[x]:x\in G\}$ is a partition of $G$.

\begin{figure}[hptb]
  \centering
  \includegraphics[width=10cm]{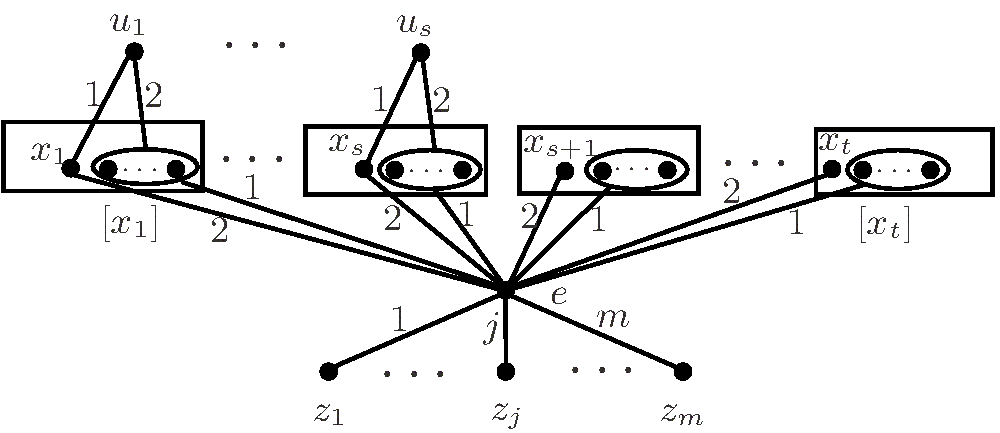}\\
  \caption{The set of edges $E_1\cup E_2\cup\{\{e,z_j\}:j=1,\ldots,m\}$}\label{e1e2}
\end{figure}

\begin{lem}\label{up}
  ${\rm rc}(\Gamma_{G})\le \max\{|M_G|,3\}$.
\end{lem}
\proof Suppose that
$\{[x_1],\ldots,[x_s]\}$ and $\{[x_{s+1}],\ldots,[x_{s+t}]\}$
are  partitions of $\{x\in G: |x|\text{ is even, } |x|\ge 4\}$
and $\{x\in G: |x|\text{ is odd, } |x|\ge 3\}$, respectively.
For $1\le i\le s$, let $u_i$ be the involution in $\langle x_i\rangle$.
Write $M_G=\{z_1,\ldots,z_m\}$ and
\begin{eqnarray*}
E_1&=&\{\{e,x\}: x\in \cup_{i=1}^{s+t}([x_i]\setminus \{x_i\})\} \cup
\{\{u_i,x_i\}: i=1,\ldots,s\},\\
E_2&=&\{\{e,x_i\}: i=1,\ldots,s+t\} \cup
(\cup_{i=1}^s\{\{u_i,x\}:x\in [x_i]\setminus \{x_i\}\}),\\
E_3&=&E(\Gamma_G)\setminus(E_1 \cup E_2\cup\{\{e,z_j\}: j=1,\ldots,m\}).
\end{eqnarray*}
The sets of edges $E_1,E_2$ and $\{\{e,z_j\}: j=1,\ldots,m\}$ are showed in Figure~\ref{e1e2}.

Let $k=\max\{|M_G|,3\}$. Define a coloring
\begin{equation*}\label{zeta}
\zeta: E(\Gamma_G)\longrightarrow\{1,\ldots,k\},\quad f\longmapsto
\left\{
  \begin{array}{lll}
    i, & \hbox{if } f\in E_i,&\hbox{where }i=1,2,3; \\
    j, & \hbox{if } f=\{e,z_j\},&\hbox{where }j=1,\ldots,m.
  \end{array}
\right.
\end{equation*}
In order to get the desired inequality, we only need to show that $\zeta$ is a rainbow $k$-coloring of $\Gamma_G$. Pick a pair of nonadjacent vertices $v$ and $w$ of $\Gamma_G$. It suffices to find a rainbow path from $v$ to $w$ under the coloring $\zeta$.
If $\zeta(\{e,v\})\neq\zeta(\{e,w\})$, then $(v,e,w)$ is a desired rainbow path. Now suppose $\zeta(\{e,v\})=\zeta(\{e,w\})$. Then $\{v,w\}\not\subseteq (M_G\cup\{e\})$.
Without loss of generality, we may assume that $v\in V(\Gamma_G)\setminus (M_G\cup\{e\})$.
As shown in Figure~\ref{e1e2},
there exists a vertex $v'\in V(\Gamma_G)\setminus (M_G\cup\{e\})$ such that
$$
\{\zeta(\{e,v\}),\zeta(\{e,v'\}),\zeta(\{v,v'\})\}=\{1,2,3\},
$$
which implies that $(v,v',e,w)$ is  a rainbow path, as desired.
$\qed$

Combining Lemmas~\ref{low} and~\ref{up}, we get the following.

\begin{prop}\label{mge3}
 If $|M_G|\ge 3$, then ${\rm rc}(\Gamma_{G})= |M_G|$.
\end{prop}

For a prime $p$, let $s_p(G)$ denote the number of subgroups of order $p$ in $G$.

\begin{lem}{\rm (\cite[Section 4, I]{Fr})}\label{basic}
Let  $p$ be a prime dividing the order of $G$. Then
$$s_{p}(G)\equiv1\pmod{p}.$$
\end{lem}

\begin{lem}\label{claim}
 Let $p$ be a prime dividing $|G|$. If ${\rm rc}(\Gamma_G)=2$, then  $s_p(G)=1$.
\end{lem}
\proof
Suppose for the contrary that $s_p(G)\neq 1$. It follows from Lemma~\ref{basic} that $s_p(G)\geq 3$.
Let $\langle y_1\rangle$, $\langle y_2\rangle$ and $\langle y_3\rangle$ be pairwise distinct subgroups of order $p$ in $G$. Note that, for $i\neq j$, there is no cyclic subgroup containing $\langle y_i\rangle$ and $\langle y_j\rangle$. Hence, the path from $y_i$ to $y_j$ with length $2$ is unique, which is $(y_i,e,y_j)$. For any  rainbow $k$-coloring $\zeta$ of $\Gamma_G$, we deduce that $\zeta(\{e,y_1\}),\zeta(\{e,y_2\})$ and $\zeta(\{e,y_3\})$ are pairwise distinct, which implies that $k\geq 3$,  contrary to ${\rm rc}(\Gamma_G)=2$.
$\qed$

By Lemmas~\ref{low},~\ref{up} and~\ref{claim}, we get the following result.

\begin{prop}\label{m=2}
If $|M_G|=2$, then ${\rm rc}(\Gamma_{G})= 3$.
\end{prop}

\begin{prop}\label{m=1}
 If $|G|\ge 3$ and $|M_G|=1$, then  ${\rm rc}(\Gamma_{G})= 3$.
\end{prop}
\proof
It follows from Lemma~\ref{up} that ${\rm rc}(\Gamma_{G})\le 3$. Suppose for the contrary that ${\rm rc}(\Gamma_{G})\le 2$. If ${\rm rc}(\Gamma_{G})=1$, then $\Gamma_G$ is a complete graph, and so $G$ is a cyclic  group of prime power order by \cite[Theorem 2.12]{CGS}, contrary to $|G|\ge 3$ and $|M_G|=1$. In the following assume that ${\rm rc}(\Gamma_{G})=2$.

Suppose that $G$ is a $2$-group. By Lemma~\ref{claim}, the involution is unique, which implies that $G$ is cyclic or  generalised quaternion by \cite[Theorem 5.4.10 (ii)]{Gor}, a contradiction.

Suppose that $|G|$ has a prime divisor $p$ at least $3$. Let $x$ be an element of $G$ with $|x|=p$. Write $M_G=\{z\}$.
It follows from Lemma \ref{claim} that $\langle x\rangle$ and $\langle z\rangle$ are normal subgroups  in $G$.
Note that $\langle x\rangle\cap \langle z\rangle=\langle e\rangle$. So $\langle x\rangle \langle z\rangle$ is a cyclic group,  contrary to the fact that $z$ is maximal.
$\qed$

Proof of Theorem \ref{thm1} follows from Propositions \ref{mge3}, \ref{m=2} and \ref{m=1}.

\medskip

For $n\ge 3$, let $D_{2n}$ denote the dihedral group of order $2n$, where $$D_{2n}=\langle a,b: a^n=b^2=1, bab=a^{-1}\rangle.$$
$\mathbb{Z}_2^n$ denotes the elementary abelian $2$-group.  Note that $M_{D_{2n}}=\{b,ab,a^2 b,\ldots,a^{n-1}b\}$ and $M_{\mathbb{Z}_2^n}$ consists of all nonidentity elements.   By Theorem~\ref{thm1}, we get the following.
\begin{eg}\label{D2n}
For $n\ge 3$, we have ${\rm rc}(\Gamma_{D_{2n}})=n$ and $\rc(\Gamma_{\mathbb{Z}_2^n})=2^n-1$.
\end{eg}

\section{$M_G=\emptyset$}\label{MG=0}

 In this section we study the rainbow connection number of $\Gamma_G$ when $G$ has no maximal involutions.

\begin{figure}[hptb]
  \centering
  \includegraphics{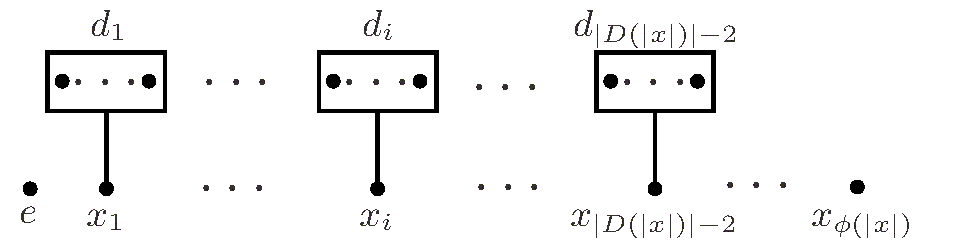}\\
  \caption{The partition of $V(\Gamma_{\langle x\rangle})$ and the set of edges $E_1(\langle x\rangle)$}\label{e1x}
\end{figure}

 For a positive integer $n$, let $D(n)$ be the set of all divisors of $n$. Denote by $\phi$ the Euler's totient function. In view of \cite[Part \uppercase\expandafter{\romannumeral8}, Problem 45]{PS76}, one has $\phi(n)\ge |D(n)|-2$.
 For $x\in G$, recall that $[x]=\{y\in G:\langle y\rangle=\langle x\rangle\}$. Write
\begin{equation}\label{1}
E_1(\langle x\rangle)=\bigcup_{i=1}^{|D(|x|)|-2}\{\{x_i,y\}:y\in\langle x\rangle,|y|=d_i\},
\end{equation}
where $[x]=\{x_1,\ldots,x_{\phi(|x|)}\}$ and $D(|x|)=\{1, d_1,\ldots,d_{|D(|x|)|-2},|x|\}$ (see Figure~\ref{e1x}).
\begin{thm}\label{thm2}
  Let $G$ be a finite group with no maximal involutions.

  {\rm(i)} If $G$ is cyclic, then
  $$
\rc(\Gamma_{G})=\left\{
                                  \begin{array}{ll}
                                    1, & \hbox{if $|G|$ is a prime power;} \\
                                    2, & \hbox{otherwise.}
                                  \end{array}
                                \right.
$$

{\rm(ii)} If $G$ is noncyclic, then $\rc(\Gamma_G)=2$ or $3$.
\end{thm}
\proof (i) Write $G=\langle x\rangle$. If $|x|$ is a prime power, then $\Gamma_{G}$ is a complete graph by \cite[Theorem 2.12]{CGS}, and so ${\rm rc}(\Gamma_{G})=1$. Now suppose that $|x|$ is not a power of any prime. Then ${\rm rc}(\Gamma_{G})\ge 2$.
With reference to (\ref{1}), write $E_1=E_1(\langle x\rangle)$.
It is clear that $E_1\subseteq E(\Gamma_{G})$. Let $E_2=E(\Gamma_{G})\setminus E_1$. Define a coloring
$$
\zeta: E(\Gamma_{G})\longrightarrow\{1,2\},\quad f\longmapsto i\text{ if }f\in E_i.
$$
In order to get the desired result, we only need to show that $\zeta$ is a rainbow $2$-coloring.
For any pair of nonadjacent vertices $v$ and $w$, there exist distinct indices $i$ and $j$ in $\{1,\ldots,|D(|x|)|-2\}$ such that $|v|=d_i$ and $|w|= d_j$. It follows from Figure~\ref{e1x} that $(v,x_i,w)$ is a rainbow path under the coloring $\zeta$, as desired.

(ii) It is immediate from Lemma~\ref{up}.
$\qed$

We first give two examples for computing $\rc(\Gamma_G)$ when $G$ is noncyclic with no maximal involutions.
The generalized quaternion group  is defined by
\begin{equation}\label{2}
Q_{4n}=\langle x,y: x^n=y^2, x^{2n}=1, y^{-1}xy=x^{-1}\rangle,\qquad n\ge 2.
\end{equation}
\begin{eg}\label{rq8zn}
  If $n$ is odd, then $\rc(\Gamma_{Q_8\times \mathbb{Z}_n})=2$.
\end{eg}
\proof There are exactly three  maximal cyclic subgroup in $Q_8\times \mathbb{Z}_n$, which we denote by $\langle x_1\rangle$, $\langle x_2\rangle$ and $\langle x_3\rangle$. It is easy to see that $|x_1|=|x_2|=|x_3|=4n$. Let $C$ be a subgroup of order $2n$ in $\langle x_1\rangle$. Then $C=\langle x_i\rangle\cap \langle x_j\rangle$ for $1\le i<j\le 3$. Write $D(n)=\{d_1,\ldots,d_t\}$. Let $B_i$, $C_i$ and $D_i$ be the set of generators of the subgroup of order $4d_i$ in $\langle x_1\rangle$, $\langle x_2\rangle$ and $\langle x_3\rangle$, respectively.   Consequently, we have
\begin{eqnarray*}
  V(\Gamma_{Q_8\times \mathbb{Z}_n})&=&C\cup\bigcup_{i=1}^t(B_i\cup C_i\cup D_i),\\
  E(\Gamma_{Q_8\times \mathbb{Z}_n})&=& E(\Gamma_{\langle x_1\rangle})\cup E(\Gamma_{\langle x_2\rangle})\cup E(\Gamma_{\langle x_3\rangle}).
\end{eqnarray*}
\begin{figure}[hptb]
  \centering
  \includegraphics{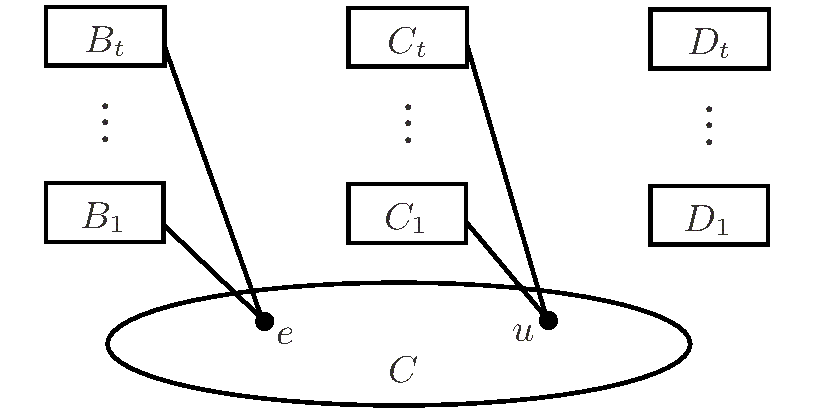}\\
  \caption{The partition of $V(\Gamma_{Q_8\times \mathbb{Z}_n})$ and the set of edges  $E_1'$}\label{q8zn}
\end{figure}
The partition of $V(\Gamma_{Q_8\times \mathbb{Z}_n})$ is showed in Figure~\ref{q8zn}, where $u$ is the unique involution.  With reference to (\ref{1}), there exists a unique vertex $x_3'\in[x_3]$ such that $\{u,x_3'\}\in E_1(\langle x_3\rangle)$.
Write
\begin{eqnarray*}
E_1'&=&\bigcup_{i=1}^t(\{\{e,x\}:x\in B_i\}\cup\{\{u,x\}:x\in C_i\}),\\
E_1&=&E_1'\cup E_1(\langle x_1\rangle)\cup E_1(\langle x_2\rangle)\cup (E_1(\langle x_3\rangle)\setminus\{\{u,x_3'\}\}).
\end{eqnarray*}
It is clear that $E_1\subseteq E(\Gamma_{Q_8\times \mathbb{Z}_n})$. Write $E_2=E(\Gamma_{Q_8\times \mathbb{Z}_n})\setminus E_1$. Define a coloring
$$
\zeta: E(\Gamma_{Q_8\times \mathbb{Z}_n})\longrightarrow\{1,2\},\quad f\longmapsto k\text{ if }f\in E_k.
$$

For $i=1,2,3$, let $\Delta_i$ be the subgraph of $\Gamma_{\langle x_i\rangle}$ induced by $V(\Gamma_{\langle x_i\rangle})\setminus\{e,u\}$.  Similar to the proof of Theorem~\ref{thm2} (i), we deduce that $\zeta|_{E(\Delta_i)}$ is a rainbow $2$-coloring of $\Delta_i$. If vertices $v$ and $w$ satisfy $u\not\in\{v,w\}$ and  $\{v,w\}\not\subseteq V(\Delta_i)$ for any $i\in\{1,2,3\}$, then $(v,e,w)$ or $(v,u,w)$ is a rainbow path under $\zeta$ from Figure~\ref{q8zn}. If $v$ is a vertex that  is not adjacent to $u$, there exists a vertex $x_3''\in [x_3]\setminus\{x_3'\}$ such that $\{x_3'',v\}\in E_1(\langle x_3\rangle)$, and so  $(u,x_3'',v)$ is a rainbow path under $\zeta$.
It follows that $\zeta$ is a rainbow $2$-coloring of $\Gamma_{Q_8\times\mathbb Z_n}$.
This completes the proof.
$\qed$

\begin{eg}\label{q4n}
If $n\ge 3$, then $\rc(\Gamma_{Q_{4n}})=3$.
\end{eg}
\proof With reference to (\ref{2}), we have
$y^{-1}=x^ny$ and $(x^iy)^{-1}=x^{2n-i}y$ for $i\in\{1,\ldots,n-1\}$, which implies that
\begin{eqnarray*}
V(\Gamma_{Q_{4n}})&=& \{e,x,\ldots,x^{2n-1}\}\cup(\bigcup_{i=0}^{n-1}\{x^iy,(x^iy)^{-1}\}),\\
E(\Gamma_{Q_{4n}})&=&E(\Gamma_{\langle x\rangle})\cup\bigcup_{i=0}^{n-1}E(\Gamma_{\langle x^iy\rangle}),
\end{eqnarray*}
\begin{figure}[hptb]
  \centering
  \includegraphics{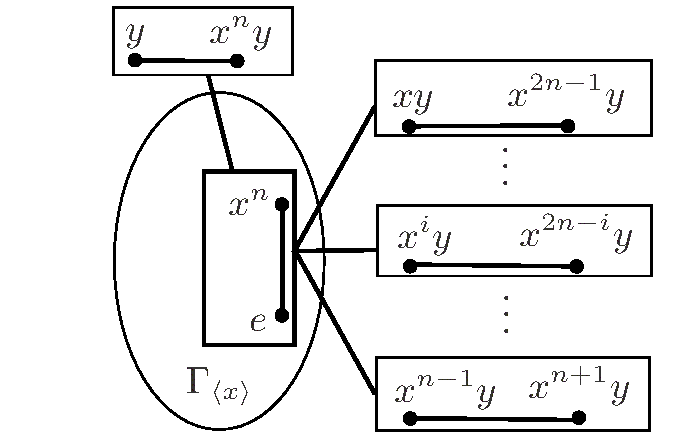}\\
  \caption{$\Gamma_{Q_{4n}}$}\label{gq4n}
\end{figure}
as shown in Figure~\ref{gq4n}. It follows from Theorem~\ref{thm2} that $\rc(\Gamma_{Q_{4n}})=2$ or $3$.
Suppose  for the contrary that there exists a rainbow $2$-coloring $\zeta$ of $\Gamma_{Q_{4n}}$.

Assume that $n=3$. Without loss of generality, let $\zeta(\{e,x^2\})=1$. Then $\zeta(\{e,x^iy\})=2$ for $i\in\{0,1,2\}$. Hence, for $0\le i<j\le 2$, the rainbow path from $x^iy$ to $x^jy$ is $(x^iy,x^3,x^jy)$, which implies that $\zeta(\{y,x^3\})$, $\zeta(\{xy,x^3\})$ and $\zeta(\{x^2y,x^3\})$ are pairwise distinct, a contradiction. Therefore $\rc(\Gamma_{Q_{12}})=3$.

In the following, assume that $n\ge 4$. Let $\Delta$ be the induced subgraph of $\Gamma_{Q_{4n}}$ on the vertices $\{e,x,y,xy,x^2y,x^3y,x^n\}$.  Then $\zeta|_{E(\Delta)}$ is  a rainbow $2$-coloring of $\Delta$.

We claim that there exists a rainbow path from $e$ to $x^n$ with length $2$ under $\zeta|_{E(\Delta)}$ in $\Delta$.
In fact, if $\zeta|_{E(\Delta)}(\{e,x^iy\})=\zeta|_{E(\Delta)}(\{x^iy,x^n\})$ for each $i\in\{0,1,2,3\}$, then there exist two distinct indices $j$ and $k$ in $\{0,1,2,3\}$ such that
$$\zeta|_{E(\Delta)}(\{e,x^jy\})=\zeta|_{E(\Delta)}(\{x^jy,x^n\})
=\zeta|_{E(\Delta)}(\{e,x^ky\})=\zeta|_{E(\Delta)}(\{x^ky,x^n\}),$$
which implies that there is no rainbow path from $x^jy$ to $x^ky$ under
$\zeta|_{E(\Delta)}$ in $\Delta$, a contradiction. Hence, the claim is valid.

Let $\Delta_0$ be the graph obtained from $\Delta$ by deleting the edge $\{e,x_n\}$. Then $\Delta_0$ is isomorphic to the complete bipartite graph $K_{2,5}$. By the claim above, we have $\rc(K_{2,5})=2$, contrary to \cite[Theorem 2.6]{CJMZ}.
$\qed$

For a noncyclic group $G$ with no maximal involutions, it is difficult for us to determine which groups $G$ satisfy $\rc(\Gamma_G)=2$. However, we give a sufficient condition.

\begin{figure}[htpb]
\begin{center}
  \includegraphics[width=10cm]{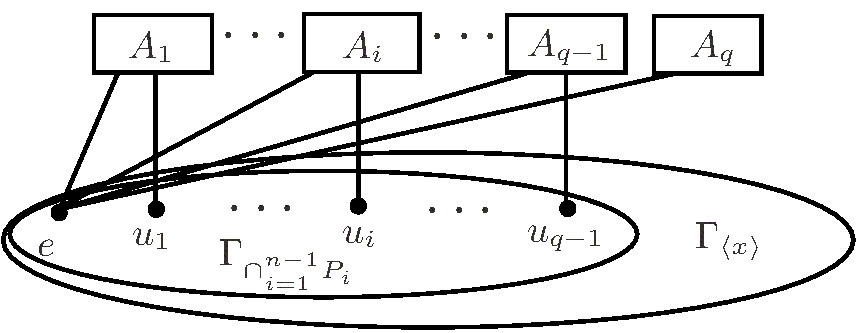}\\
  \caption{ $V(\Gamma_G)$ and the set of edges $E_1'$}\label{pnq}
\end{center}
\end{figure}

\begin{prop}\label{suf2}
If $G$ is a group of order $p^nq$ for positive integer $n$, where $p,q$ are distinct primes and $p<q$,  such that the following conditions hold,  then $\rc(\Gamma_G)=2$.

{\rm(i)} Each Sylow $p$-subgroup is cyclic and the Sylow $q$-subgroup is unique.

{\rm(ii)} The intersection of all Sylow $p$-subgroups is of order $p^{n-1}$.

{\rm(iii)} $p^{n-1}\ge q$.
\end{prop}
\proof Note that the number of Sylow $p$-subgroups is $q$. Suppose that $\{P_1,\ldots, P_q\}$ is the set of all Sylow $p$-subgroups, and $Q$ is the unique Sylow $q$-subgroup. Then $\cap_{i=1}^qP_i$ and $Q$ are cyclic and normal in $G$. Hence, there exists an element $x$ of order $p^{n-1}q$ such that  $(\cap_{i=1}^qP_i)Q=\langle x\rangle$,  and so the set of all cyclic subgroups of $G$ is
$$
\{P_1,\ldots,P_q\}\cup\{\langle y\rangle: y\in\langle x\rangle\}.
$$
For $1\le i\le q$,  let $A_i$ be the set of all generators of $P_i$.  By (iii) we choose   pairwise distinct elements $u_1,\ldots, u_{q-1}$ in $(\cap_{i=1}^qP_i)\setminus\{e\}$.
With reference to (\ref{1}), write
\begin{eqnarray*}
  E_1'&=&\{\{e,y\}: y\in \bigcup_{i=1}^qA_i\}\cup\bigcup_{i=1}^{q-1}\{\{u_i,y\}: y\in A_i\},\\
   E_1&=&E_1'\cup E_1(\langle x\rangle).
\end{eqnarray*}
The set $E_1'$ is showed in Figure~\ref{pnq}.
It is clear that $E_1\subseteq E(\Gamma_G)$. Let $E_2=E(\Gamma_{G})\setminus E_1$. Define a coloring
$$
\zeta: E(\Gamma_{G})\longrightarrow\{1,2\},\quad f\longmapsto k\text{ if }f\in E_k.
$$

In order to get the desired result, we only need to show that $\zeta$ is a rainbow $2$-coloring of $\Gamma_G$. It follows from Theorem~\ref{thm2}   that $\zeta|_{E(\Gamma_{\langle x\rangle})}$ is a rainbow $2$-coloring of $\Gamma_{\langle x\rangle}$.
Pick any pair of nonadjacent vertices $z$ and $w$ such that $\{z,w\}\not\subseteq V(\Gamma_{\langle x\rangle})$.  It suffices to find a rainbow path from $z$ to $w$ under $\zeta$. Without loss of generality, assume that $z\in \cup_{i=1}^qA_i$.
If $w\in\cup_{i=1}^qA_i$, then there exist indices $i$ and $j$ in $\{1,\ldots,q\}$ with $i<j$ such that $z\in A_i$ and $w\in A_j$, and so $(z,u_i,w)$ is a desired rainbow path. If $w\in V(\Gamma_{\langle x\rangle})$,  then $(z,e,w)$ is a desired rainbow path.
$\qed$

By Proposition~\ref{suf2}, we have the following example.

\begin{eg}
  Let $G=\langle a,b:a^{27}=b^7=e, a^{-1}ba=b^2\rangle\cong\mathbb{Z}_{27}\ltimes \mathbb{Z}_7$. Then $\rc(\Gamma_G)=2$.
\end{eg}

The following sufficient condition for $\rc(\Gamma_G)=3$ is immediate from Theorem~\ref{thm2} and   Lemma~\ref{claim}.
\begin{prop}\label{suf}
 Suppose that $G$ is a noncyclic group with no maximal involutions.
 If  there exists a prime $p$ dividing $|G|$ such that the subgroup of order $p$ in $G$ is not unique, then ${\rm rc}(\Gamma_{G})= 3$.
\end{prop}

Finally, we determine the rainbow connection number of the power graph of a nilpotent group.

\begin{cor}\label{}
Let $G$ be a noncyclic nilpotent group with no maximal involutions. Then
$$
{\rm rc}(\Gamma_G)=\left\{
  \begin{array}{ll}
    2, & \hbox{if $G$ is isomorphic to  $Q_8\times \mathbb{Z}_n$ for some odd number $n$;} \\
    3, & \hbox{otherwise.}
  \end{array}
\right.
$$
\end{cor}
\proof It follows from Theorem~\ref{thm2} that $\rc(\Gamma_G)=2$ or $3$.  Suppose $\rc(\Gamma_G)=2$. Then for any prime $p$ dividing $|G|$, the subgroup of order $p$ in $G$ is unique by Proposition~\ref{suf}. By \cite[Theorem 5.4.10 (ii)]{Gor}, the Sylow $p$-subgroups are cyclic for any odd prime $p$, which implies that $2$ is a divisor of $|G|$ and the Sylow $2$-subgroup is isomorphic to $Q_{2^m}$ for $m\geq 3$. Hence we get $G\cong Q_{2^m}\times \mathbb{Z}_n$ for some odd number $n$. Let $H$ be a subgroup of $G$ that is isomorphic to $Q_{2^m}$.

We claim that for any pair of nonadjacent vertices $x$ and $y$ of $\Gamma_H$, there does not exist a vertex  in $G\setminus H$ adjacent to both $x$ and $y$ in $\Gamma_G$. Suppose for the contrary that $\{\{x,z\},\{y,z\}\}\subseteq E(\Gamma_G)$ for some $z\in G\setminus H$. Then $x=z^s$ and $y=z^t$ for some integers $s$ and $t$, which implies that $x,y\in\langle z\rangle$. Note that $|x|$ and $|y|$ are powers of $2$. It follows that
$|x|$ is divisible by $|y|$, or $|y|$ is divisible by $|x|$.
Since $x,y\in\langle z\rangle$, one has $\langle x\rangle \subseteq \langle y\rangle$ or $\langle y\rangle \subseteq \langle x\rangle$. Thus  $x$ and $y$ are adjacent, a contradiction. Hence, the claim is valid.
By the claim, one gets $\rc(\Gamma_H)=2$. It follows from Example~\ref{q4n} that $m=3$, and so $G\cong Q_8\times \mathbb{Z}_n$. By Example~\ref{rq8zn}, we get the desired result.
$\qed$

\section*{Acknowledgement}
The authors are grateful to the referees for many useful suggestions and comments.
This research is supported by National Natural Science Foundation of China (11271047,
11371204) and the Fundamental Research Funds for the Central University of China.

\end{CJK*}

\end{document}